\documentclass[reqno]{amsart}

\usepackage{amsmath,amsthm,amssymb,amsfonts,mathrsfs,amscd}
\usepackage{multirow}
\usepackage{graphicx}

\usepackage{bm}
\usepackage{mathrsfs}
\usepackage{amsmath}
\usepackage{graphicx}
\usepackage{amssymb}
\usepackage{dsfont}
\usepackage{tikz}
\usepackage{graphicx}
\usepackage{tikz-3dplot}
\usepackage{subfig}

\newtheorem{remark}{Remark}[section]

\newcommand{\stl}
{\mathrel{\raise2pt\hbox{${\mathop<\limits_{\raise1pt\hbox
{\mbox{$\sim$}}}}$}}}

\newcommand{\ste}
{\mathrel{\raise2pt\hbox{${\mathop=\limits_{\raise1pt\hbox{\mbox{$\sim$}}}}$}}}

\newcommand{\stg}
{\mathrel{\raise2pt\hbox{${\mathop>\limits_{\raise1pt\hbox{\mbox{$\sim$}}}}$}}}

\newcommand{\ms}[1]{\mathscr{#1}}

\makeatletter
  \newcommand\figcaption{\def\@captype{figure}\caption}
  \newcommand\tabcaption{\def\@captype{table}\caption}
\makeatletter

\def\su{\sum\limits^{N_0}_{i=1}}

\def\be{\begin{equation}}
\def\bn{\begin{eqnarray}{}}
\def\ee{\end{equation}}
\def\en{\end{eqnarray}}
\def\bq{\begin{eqnarray*}{}}
\def\eq{\end{eqnarray*}}
\def\l{\label}
\def\ri{\rightarrow}
\def\o{\omega}
\def\O{\Omega}
\def\la{\lambda}
\def\al{\alpha}

\def\T{{\mathcal T}_h}

\def\G{\Gamma}

\def\p{\partial}

\def\12{{1\over 2}}

\def\a{{\mathcal A}}
\def\q{\quad}

\def\V0{\bar{V}_h}

\def\ff{ {\mbox{\sc f}}}
\def\se{{\mbox{\sc e}}}
\def\vv{\mbox{\sc v}}
\def\o{\Omega}
\def\c{{\bf curl}}
\def\u{{\bf u}}
\def\v{{\bf v}}

\def\x{{\bf x}}
\def\ti{\times}

\def\w{{\bf w}}
\def\n{{\bf n}}
\def\f{{\bf f}}
\def\0{{\bf 0}}
\def\a{{\bf a}}
\def\b{{\bf b}}
\def\t{{\bf t}}
\def\r{{\bf r}}

\def\12{{1\over 2}}

\def\su{\sum\limits}

\def\rig{\rightarrow}

\def\chi{{\mathcal X}}

\def\A{{\mathcal A}}

\begin{document}

\title{Substructuring preconditioners with novel interface solvers for
general elliptic-type equations in three dimensions}

\author{Qiya Hu}
\author{Shaoliang Hu}

\thanks{LSEC, Institute of Computational
Mathematics and Scientific Engineering Computing,  Chinese Academy
of Sciences, Beijing 100080, China (hqy@lsec.cc.ac.cn) \and hushaoliang@lsec.cc.ac.cn).
This work was funded by Natural Science Foundation of China G11571352.}

\maketitle
\begin{abstract}
In this paper we propose two variants of the substructuring preconditioner for solving
three-dimensional elliptic-type equations with strongly discontinuous coefficients.
In the new preconditioners, we use the simplest coarse solver associated with the finite element space induced by the coarse
partition, and construct novel interface solvers based on some new observations. The resulting preconditioners share
the merits of the non-overlapping domain decomposition method (DDM) and the overlapping DDM in the sense that they not only are cheap but also are easy to implement.
We apply the proposed preconditioners to solve the linear elasticity problems and Maxwell's equations in three dimensions.
Numerical results show that the convergence rate of PCG method
with the preconditioners are nearly optimal, and also robust with respect to
the (possibly large) jumps of the coefficients in the considered equations.
\end{abstract}



\pagestyle{myheadings}
\thispagestyle{plain}
\markboth{Substructuring preconditioners  for elliptic-type equations}{QIYA HU AND SHAOLIANG HU}

\section{Introduction}
There are many works to study (non-overlapping or overlapping) domain decomposition methods (DDMs) for solving the systems generated by finite element discretization of
elliptic-type partial differential equations (\cite{BrambleP1989}-\cite{CaiS1999},\cite{ChanZ1994}-\cite{Gander2012}, \cite{Ha1991}, \cite{Hu2004}-\cite{Mandel2005}, \cite{Smith1992}-\cite{Toselli2006}, \cite{VeigaCPS2012}-\cite{VeigaPSWZ2014},\cite{re25,XuZ1998} and the references therein). Non-overlapping DDMs and overlapping DDMs have their respective merits and drawbacks: non-overlapping DDMs are cheaper and more efficient for the case of large jump coefficient than overlapping DDMs (with large overlap), but non-overlapping DDMs are more difficult to construct and implement than overlapping DDMs. In fact, the construction of non-overlapping DDMs heavily depends on the considered models. For example, non-overlapping DDMs for positive definite Maxwell's equations are essentially different from that for the usual elliptic equation (comparing  \cite{DohrmannW2015,HuSZ2011,HuZ2004,Toselli2006}).
The drawbacks mentioned above restrict applications of the non-overlapping DDMs and the overlapping DDMs with large overlap. Although the overlapping DDMs with small overlap are cheap and easy to implement, they have slower
convergence. Over the past two decades, some interesting DDMs have been proposed and analyzed, for example, the DDMs with Lagrangian multipliers \cite{Farhat1991, Farhat2000, HuS2010, KlawWD2002}
, the BDDC methods \cite{Dohrmann2003, DohrmannW2015, LiW2007, Mandel2003}, the restricted additive Schwarz methods \cite{CaiS1999,FrommerS2001}, the optimized Schwarz methods \cite{Gander2006,Gander2012}. These methods have obvious advantages over the traditional DDMs: the DDMs with Lagrangian multipliers can be conveniently handle non-matching grids, the BDDC methods are particularly practical for the case with irregular subdomains, the restricted additive Schwarz methods are cheaper and faster than the standard overlapping DDMs, the optimized Schwarz methods can accelerate convergence of the non-overlapping Schwarz methods.

In the present paper, we try to construct relatively united substructuring preconditioners for elliptic-type equations, such that they are cheap, easy to implement and have fast convergence.
As usual, we decompose the considered domain into the union of some non-overlapping subdomains, which constitute a coarse partition of the domain.
In the proposed preconditioners, we use the simplest coarse space induced by the coarse
partition as in the overlapping DDMs. The main goal of this paper is to design cheap and practical local interface solvers based on some new observations.

For each internal cross-point, we introduce an auxiliary subdomain that contains the internal cross-point as its ``center" and has almost the same size with the original subdomains. Associated
with each auxiliary subdomain, we define a local interface problem such that the solution of the local interface problem is discrete harmonic in the intersection of the auxiliary subdomain with every original subdomain
adjoining it. Notice that each intersection is only a part of some original subdomain, so the local interface problem is defined on a space consisting of ``inexact" harmonic extensions. The corresponding local interface solver is implemented by
solving a Dirichlet problem (residual equation), which is defined on the natural restriction space of the original finite element space on the auxiliary subdomain. It is clear that each local interface solver
has almost the same cost with an original subdomain solver. We would like to point out that the proposed local interface solvers are different from the existing local interface solvers defined in the vertex space method \cite{Smith1992} or the interface overlapping additive Schwarz  \cite{XuZ1998}, where exact harmonic extensions are required.

In order to further reduce the cost of the local interface solvers described above, we need to decompose each local interface problem into two subproblems and present approximate local interface solvers based on a coarsening technique. In the step for solving a local
interface problem, we are interested only in the degrees of freedom on the local interface, instead of the degrees of freedom in the interiors of subdomains. Intuitively, the accuracy of the degrees of freedom on
the local interface are not sensitive to the grids far from the local interface. Based on this observation, we construct auxiliary non-uniform
grids in each subdomain adjoining the considered local interface such that the auxiliary grids coincide with the original fine grids on the local interface but gradually become coarser when nodes are far from the local interface.     These auxiliary grids can be easily generated by the existing software and contain much smaller number of nodes than the original fine grids in a subdomain. Such an approximate local interface solver is implemented by solving
a Dirichlet problem on the finite element space defined by the auxiliary grids, and so it is much cheaper than the original local interface solver.

The constructions of the coarse solver and the proposed local interface solvers
do not depend on the considered models, and the resulting substructuring preconditioners are cheap and easy to implement.
As pointed out in \cite{DohrmannW2015}, the design of an efficient substructuring preconditioner for three dimensional Maxwell's equations poses quite significant
challenges. A few existing preconditioners on this topic are either expensive or difficult to implement. We will apply the proposed substructuring preconditioners to solve the linear elasticity problems
and Maxwell's equations in three dimensions. Numerical results show that the preconditioners are robust uniformly for the two kinds of equations even if the coefficients have large jumps.

The outline of the paper is as follows. In Section 2, we give the variational formula of general elliptic-type equations and introduce a partition based on domain
decomposition.  In Section 3, we describe local interface solvers associated with vertex-related subdomains and define the resulting substructuring preconditioner for the general elliptic system.
In Section 4, we design cheaper local interface solvers and present the corresponding preconditioner based on a coarsening technique. In Section 5, we discuss applications of the substructuring methods
to elasticity problems and Maxwell's equations.
In section 6, we will report some numerical results for the linear elasticity problems and Maxwell's equations.

\section{Elliptic-type equations and domain decomposition}
In this section, we describe the considered problems.
\subsection{Elliptic-type equations}
Let $\Omega$ be a bounded and connected Lipschitz domain in $\mathds{R}^3$. For convenience, we just consider the weak form of elliptic-type equations. Let $V(\Omega)$ denote a Hilbert
space with the scalar product $(\cdot,\cdot)_V$, and $||\cdot||_V$ be the induced norm. We introduce a real bilinear form $\bm{\A}(\cdot,\cdot): V(\Omega) \times V(\Omega) \rightarrow R$. We assume that $\bm{\A}(\cdot,\cdot)$ is symmetric, continuous and coercive in the sense that
\[ \bm{\A}(\u,\v) = \bm{\A}(\v,\u), \quad  |\bm{\A}(\u,\v)| \leq c_1||\u||_V ||\v||_V,\quad\forall\u,\v \in V(\Omega) \]
and
\[ \bm{\A}(\u,\u) \geq c_2 ||\u||^2_V, \quad \forall\u\in V(\Omega) \]
for two positive number $c_1$ and $c_2$.

Given a linear functional $ \bm{F} \in V'(\Omega)$, we consider the following problem:
\begin{equation}
\begin{cases}
 Find \q \u \in V(\o) \quad .st. \\
\bm{\A}(\u, \v)=\langle \bm{F},\v \rangle, \quad \forall\v\,\in V(\Omega)
 \end{cases}
 \label{eq:2.1}
\end{equation}

\subsection{Domain decomposition and discretization}
For convenience, we assume that $\Omega$ is a polyhedra. For a number $d\in
(0,~1)$, let $\Omega$ be decomposed into the union of
non-overlapping tetrahedra (or hexahedra) $\{\O_k\}$ with the size
$d$. Then we get a non-overlapping domain decomposition for $\O$:
$\bar{\O}=\bigcup\limits_{k=1}^{N} \bar{\O}_k$. Assume that
$\O_i\cap\O_j=\emptyset$ when $i\not=j$; if $i\not=j$ and
$\partial\O_i\cap\partial\O_j\not=\emptyset$, then
$\partial\O_i\cap\partial\O_j$ is a common, or a common edge, or a common vertex of $\O_i$ and
 $\O_j$. It is clear that the subdomains $\O_1,\cdots,\O_N$
constitute a {\it coarse} partition ${\mathcal T}_d$ of $\O$. If
$\partial\O_i\cap\partial\O_j$ is just a common face of $\O_i$ and
 $\O_j$, then set $\G_{ij}=\partial\O_i\cap\partial\O_j$. Define $\G=\cup\G_{ij}$.
 By $\Gamma_k$ we denote the intersection of $\Gamma$ with the boundary of the subdomain
$\o_k$. So we have $\Gamma_k=\p\o_k$ if $\o_k$ is an interior
subdomain of $\o$.

With each subdomain $\O_k$ we associate a regular partition made
of tetrahedral elements (or hexahedral elements). We require that
the partitions in all the subdomains match on the common face between
two neighboring subdomains, and so they constitute a partition $\T$ on the
domain $\O$, which we assume is quasi-uniform. We denote by $h$ the mesh
 size of $\T$, i.e., $h$ denotes the maximum diameter of tetrahedra in
 the mesh $\T$.

For an element $K\in {\mathcal T}_h$, let $R(K)$ denote a set of basis functions on the element $K$. The definition
of $R(K)$ depends on the considered models, and will be given in Section 5. Define the finite element space
$$ V_{h}(\o)=\Big\{\v\in V({\Omega)}: ~ \v|_K\in R(K), ~\forall K\in {\mathcal T}_h
 \Big\}.
$$
Consider the discrete problem of (\ref{eq:2.1}):
{\it Find $\bm{u}_h\in V_{h}(\O)$ such that}
\begin{equation}
\bm{\A}(\bm{u}_h, \bm{v})=\langle \bm{F}, \bm{v} \rangle,
\quad\forall \bm{v}\in V_{h}(\O). \label{eq:2.2}
\end{equation}
This is the discrete variational problem that we need to solve in this paper.

For convenience, we define the discrete operator $A: V_h(\Omega) \rightarrow V_h(\Omega)$ as
\[ \langle A\u, \v\rangle = \bm{\A}(\u,\v), \quad \u, \v \in V_h(\Omega), \]
where $\langle \cdot, \cdot\rangle$ denotes the duality pairing between $V'(\Omega)$ and $V(\Omega)$.
Then (\ref{eq:2.2}) can be written in the operator form
\begin{equation}
Au_h=\bm{f}.
\label{eq:2.3}
\end{equation}
By the assumptions on ${\mathcal A}(\cdot,\cdot)$, the operator is symmetric and positive definite. Thus the above equation can be iteratively solved by PCG method.
In the rest of this paper, we will construct preconditioners for the operator $A$.

Before constructing the desired preconditioners, we first introduce some useful sets and subspaces.

${\mathcal N}_h$: the set of all nodes generated by the {\it fine} partition ${\mathcal T}_h$;

${\mathcal E}_h$: the set of all {\it fine} edges generated by the partition ${\mathcal
T}_{h}$;

${\mathcal F}_h$: the set of all {\it fine} faces generated by the partition
${\mathcal T}_h$;

${\mathcal N}_d$: the set of all nodes generated by the {\it coarse} partition ${\mathcal T}_d$.




In most applications, the degrees of freedom of $\v\in V_{h}(\o)$ are defined at the nodes in ${\mathcal N}_h$ (the nodal elements), or on the edges in ${\mathcal E}_h$ (Nedelec edge elements), or on the
faces in ${\mathcal F}_h$ (Raviart-Thomas face elements).  Throughout this paper, for a subset $\ff$ that is the union of faces in ${\mathcal F}_h$, the term
``the~degrees~of~freedom~of $\v$ vanish on~$\ff$" means that ``$\v$ has the zero~degrees~of~freedom~at~the nodes, or fine edges, or fine faces of $\ff$".

Let $G\subset \Omega$ be a subdomain that is the union of some elements in ${\mathcal T}_h$. Define
$$ V^0_h(G)=\{v\in V_{h}(\o):~~\mbox{the~degrees~of~freedom~of}~\v~\mbox{vanish~on}~\partial G\}.$$
For example, when $G=\Omega_k$ the space $V^0_h(\Omega_k)$ is just the subdomain space in the traditional substructuring methods.

For the construction of preconditioners, we will use the simplest coarse space $V_d(\Omega)$, which is defined as the finite element space
associated with the {\it coarse} partition ${\mathcal T}_d$ (see \cite{DryjaSW1994}, \cite{HuS2010}, \cite{HuSZ2011} and \cite{XuZ1998}). It is clear that
$V_{d}(\o)\subset V_{h}(\o)$.

\section{Preconditioner (I): with local interface solvers related to vertices}
\setcounter{equation}{0}
This section is devoted
to describing the first preconditioner, in which local interface solvers are defined in vertex-related subspaces.

\subsection{Space decomposition}
For each $\vv \in {\mathcal N}_d$, we construct an open region $\Omega_{\vv}^{half}$, whose ``center" is $\vv$ and size is about $d$. When $\vv\in\partial\Omega$,
the auxiliary subdomain
$\Omega_{\vv}^{half}$ is chosen as the part in $\Omega$. We assume that: (i) each
subdomain $\Omega_{\vv}^{half}$ is just the union of some elements in ${\mathcal T}_h$; (ii) the union of all the subdomains $\Omega_{\vv}^{half}$ is an open cover of
$\Omega$. Then all the subdomains $\Omega_{\vv}^{half}$ constitute overlapping domain decomposition of $\Omega$ (with small overlap).

\begin{figure}[ht]
\begin{center}
\includegraphics[width =7cm]{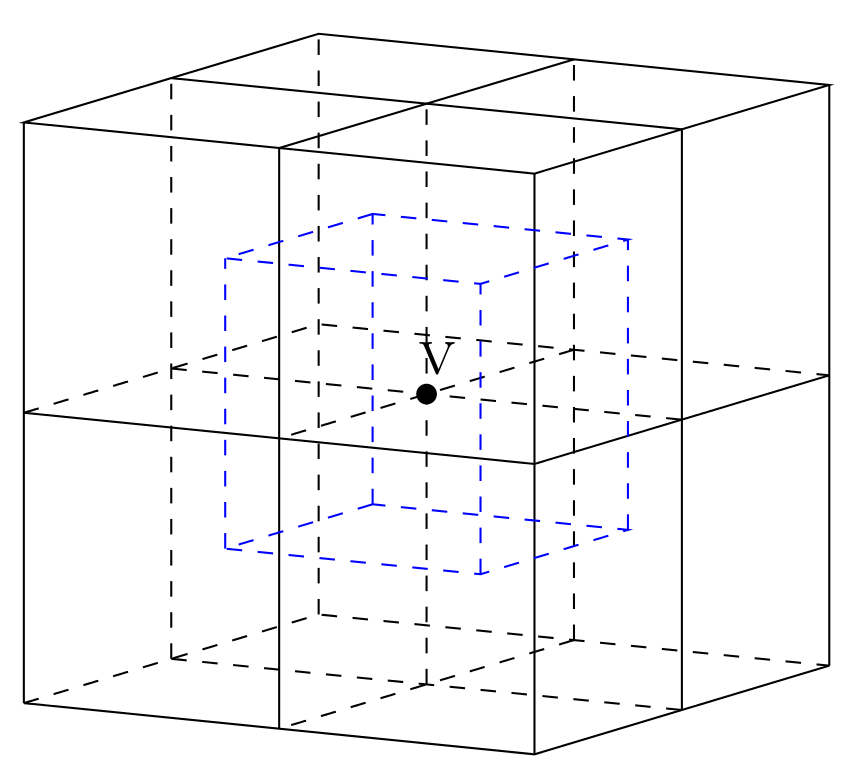}   
\caption{The auxiliary subdomain $\Omega_{\vv}^{half}$ (the blue cube) associated with the vertex $\vv$. \label{map1}}
\end{center}
\end{figure}

In order to define a space decomposition of $V_h(\Omega)$ in an exact manner, we need to introduce more notations.

For $\vv\in {\mathcal N}_d$, set
$$ \Lambda_{\vv}=\{k:~~\mbox{the~polyhedran}~\Omega_k~\mbox{contains}~\vv~\mbox{as~its~ vertex}\} $$
and define
$$ \Gamma_{\vv}^{half}=\Omega_{\vv}^{half}\cap\Gamma~~\mbox{and}~~\Omega_{\vv}=\bigcup_{k\in\Lambda_{\vv}}\Omega_k.$$
Let $V_h(\Gamma)$ denote the interface space, which consists of the natural traces of all the functions in $V_h(\Omega)$. Define the vertex-related local interface space
$$ V_h^0(\Gamma_{\vv}^{half})=\{\phi\in V_h(\Gamma):~~supp~\phi\subset \Gamma_{\vv}^{half}\}. $$
Since all the vertex-related local interfaces $\Gamma_{\vv}^{half}$ constitute an open cover of the interface $\Gamma$, we have the space decomposition
\begin{equation}
V_h(\Gamma)=\bigcup_{\vv\in{\mathcal N}_d}V^0_h(\Gamma^{half}_{\vv}).\label{3.new00}
\end{equation}

As usual, let $V_h^{\bot}(\Omega)$ denote the space consisting of all the finite element functions that are discrete $A$-harmonic in each $\Omega_k$, namely
$$ V^{\bot}_h(\Omega)=\{\v\in V_h(\o):~~{\mathcal A}(\v,\w)=0,~\forall\w\in V_h^{0}(\o_{k})~~{for~each}~~k\in {\mathcal N}_d\}.$$
Then we have
\begin{equation}
V_h(\o)=V_d(\o)+\su_{k=1}^NV_h^0(\o_k)+V^{\bot}_h(\Omega).\label{BasicDecom}
\end{equation}

For each $\vv\in {\mathcal N}_d$, define vertex-related local $A$-harmonic space
$$ V^{\bot}_h(\Omega_{\vv})=\{\v\in V^{\bot}_h(\o):~~\mbox{the~trace~of}~\v~\mbox{belongs~to}~V_h^0(\Gamma_{\vv}^{half})\}\subset V_h^0(\Omega_{\vv}).$$
In other words, $V^{\bot}_h(\Omega_{\vv})$ is just the space consisting of the discrete $A$-harmonic extensions of the functions in $V_h^0(\Gamma_{\vv}^{half})$.

It is clear that
$$ V^{\bot}_h(\Omega)=\bigcup_{\vv\in{\mathcal N}_d}V^{\bot}_h(\Omega_{\vv}). $$
Thus, by (\ref{BasicDecom}), the space $V_h(\o)$ admits the space decomposition
\begin{equation}
 V_h(\o)=V_d(\o)+\su_{k=1}^NV_h^0(\o_k)+\su_{\vv\in{\mathcal N}_d}V^{\bot}_h(\Omega_{\vv}).\label{eq:3.1}
\end{equation}

\subsection{Preconditioner}
In this subsection we define solvers on the subspaces $V_d(\o)$,
$V_h^0(\o_k)$ and $V^{\bot}_h(\o_{\vv})$.

As usual, we use $A_d: V_d(\o)\rig V_d(\o)$ and $A_k: V^0_h(\o_k)\rig V_h^0(\o_k)$ to denote the restriction of $A$ on $V_d(\o)$ and $V_h^0(\o_k)$ respectively, i.e., they satisfy
$$ (A_d\v_d,\w_d)=(A\v_d,\w_d)=\bm{\A}(\v_d,\w_d),~~~\v_d\in V_d(\o),~~\forall \w\in V_d(\o) $$
and
$$ (A_k\v, {\w})_{\o_k}=(A\v, {\w})=\bm{\A}(\v,{\w}),~~~\v\in V^0_h(\o_k),~\forall {\w}\in V_h^0(\o_k).$$

In the following we define an ``inexact" solver on $V^{\bot}_h(\o_{\vv})$. To this end, we introduce a modification of $V^{\bot}_h(\o_{\vv})$.  Let $k\in\Lambda_{\vv}$, and use $\o_{\vv,k}^{half}$ to denote
the intersection of $\o_{\vv}^{half}$ with $\Omega_k$. For each $\o_{\vv}^{half}$, define the ``inexact" $A$-harmonic space
$$ V^{\bot}_h(\o_{\vv}^{half})=\{\v\in V^0_h(\o_{\vv}^{half}):~~{\mathcal A}(\v,\w)=0,~\forall\w\in V_h^{0}
(\o_{\vv,k}^{half})~~\mbox{with}~~k\in\Lambda_{\vv}\}.$$
Notice that the functions in $V^{\bot}_h(\o_{\vv}^{half})$ have the support set $\o_{\vv}^{half}$ and are discrete $A$-harmonic only in the subdomain $\o_{\vv,k}^{half}$ of $\Omega_k$ (for any $k\in\Lambda_{\vv}$).
Thus the spaces $V^{\bot}_h(\o_{\vv}^{half})$ have essential differences from the local interface spaces proposed in the vertex space method \cite{Smith1992} or the interface overlapping additive Schwarz  \cite{XuZ1998},
where exact $A$-harmonic extensions in all $\Omega_k$ were required.

For a function $\v\in V^{\bot}_h(\o_{\vv})$, define $\v^{half}\in V^{\bot}_h(\o_{\vv}^{half})$ such that $\v^{half}=\v$ on $\Gamma_{\vv}^{half}$.
For each $\vv \in {\mathcal N}_d$, let $B_{\vv}: V^{\bot}_h(\o_{\vv})\rig V^{\bot}_h(\o_{\vv})$ be the symmetric and
positive definite operators defined by
$$
(B_{\vv}\v, \w)=\bm{\A}(\v^{half}, \w^{half}),\q\v\in V^{\bot}_h(\o_{\vv}),~~\forall\w\in V^{\bot}_h(\o_{\vv}). $$
Since the basis functions in
$V^{\bot}_h(\o_{\vv})$ are not known, the action of $B^{-1}_{\vv}$ needs to be implemented by solving a residual equation defined in $V^0_h(\o_{\vv}^{half})$
(see {\bf Algorithm 3.1} given later).

Let $Q_d: V_h(\o)\rig V_d(\o)$, $Q_k:V_h(\o)\rig V^0_h(\o_k)$
and $Q_{\vv}: V_h(\o)\rig V^{\bot}_h(\o_{\vv})$ be the standard $L^2$-projectors.
Then the first preconditioner for $A$ is defined as follows:
\be
B_{I}^{-1}=A_d^{-1}Q_d+\su_{k=1}^NA^{-1}_kQ_k+\su_{\vv\in {\mathcal N}_d}B_{\vv}^{-1}Q_{\vv}
\l{precon1} \ee
\begin{remark} To our knowledge, the coarse solver $A_d$ is the simplest and cheapest one in the non-overlapping DDMs. This coarse solver for elliptic equation was first considered in \cite{DryjaSW1994},
and then discussed in \cite{XuZ1998}. Such coarse solver was regarded as a non-optional coarse solver for long time, since the condition number of the resulting preconditioned system is not nearly optimal for the case
with large jump coefficients. Based on the framework developed in \cite{XuZ2008}, it was shown in \cite{HuS2010} that the PCG method for solving the resulting preconditioned system has the
nearly stable convergence even for the case with large jump coefficients. In \cite{HuSZ2011}, this kind of coarse solver was also applied to Maxwell's equations. When this coarse solver are used, cheap ``edge" solvers
(and ``face" solvers) need to be designed. It can be seen, from \cite{HuSZ2011} and \cite{HuS2010} (see also \cite{DryjaSW1994} and \cite{XuZ1998}), that the constructions of the existing ``edge" solvers are based on
estimates of the norms induced from the interface operators restricted on the edges and so depend on the considered models. In the proposed preconditioner $B_I$, the construction of the edge solvers $B_{\vv}$ (which also play
the role of face solvers) is unified and independent of the bilinear $\bm{\A}(\cdot,\cdot)$.
\end{remark}

\begin{remark}
The preconditioner (\ref{precon1}) can be extended to the case with irregular subdomains (i.e., $\Omega_k$ is not a polyhedron with finite faces), for which the coarse space $V_d(\Omega)$ needs to be replaced by
the image of the interpolation operator from an auxiliary regular coarse space into $V_h(\Omega)$ as in \cite{Cai1995} and \cite{ChanZ1994}.
\end{remark}

The action of the preconditioner $B^{-1}_I$, which is needed in each iteration step of
PCG method, can be described by the following algorithm.\\
{\bf Algorithm 3.1}. For ${\bf g}\in V_h(\o)$, we can compute
$\u=B^{-1}_{I}{\bf g}$
in four steps.

Step 1. Solve the system of $\u_d\in V_d(\o)$:
$$ (A_d\u_d, \v_d)=({\bf g},\v_d),~~~\forall\v_d\in V_d(\o); $$

Step 2. Solve the systems of  $\u_k\in V_h^0(\o_k)$ ($k=1,\cdots,N$) in parallel:
$$ (A_k\u_k, \v)=({\bf g},\v),~~~\forall\v\in V^0_h(\o_k),~~k=1,\cdots,N;
$$

Step 3. Solve the systems of $\u_{\vv}\in V_h^0(\o_{\vv}^{half})$ ($\vv \in {\mathcal N}_d$) in parallel:
$$ (B_{\vv}\u_{\vv}, \v)=({\bf g},\v)-\su_{k\in\Lambda_{\vv}}(A_k\u_k,\v)_{\o_k},~~~\forall\v\in
V^0_h(\o_{\vv}^{half});$$

Step 4. Compute the trace $\Phi_h=\bm{\gamma}_{\Gamma}(\su_{\vv \in {\mathcal N}_d}\u_{\vv})$, and then
compute the $A$-harmonic extension of $\Phi_h$ on each $\o_k$ to obtain
$\u^{\bot}\in V^{\bot}_h(\o)$. This leads to
$$ \u=\u_d+\su_{k=1}^N\u_k+\u^{\bot}. $$

\begin{remark} It can be seen from {\bf Algorithm 3.1} that the preconditioner (\ref{precon1}) is easy and cheap to implement (each vertex-related space $V^0_h(\o_{\vv}^{half})$ has almost the same
degrees of freedom with an original subdomain space $V^0_h(\o_k)$).

\end{remark}

\section{Preconditioner (II): with approximate interface solvers}
\setcounter{equation}{0}
Although the local interface solvers $B_{\vv}$ defined in the last section is not expensive, we want to further reduce the cost for implementing the action of $B^{-1}_{\vv}$.
To this end, we introduce a coarsening technique for the construction of cheaper local interface solvers.


${\mathcal F}_d$: the set of all the open (coarse) faces generated by the partition ${\mathcal T}_{d}$;

${\mathcal E}_d$: the set of all the open (coarse) edges generated by the partition ${\mathcal T}_{d}$;

${\mathcal F}_{\vv}$: the set of the (coarse) faces, each of which belongs to ${\mathcal F}_d$ and contains $\vv$ as its vertex;

${\mathcal E}_{\vv}$: the set of the (coarse) edges, each of which belongs to ${\mathcal E}_d$ and contains $\vv$ as its vertex;

For $\se\in {\mathcal E}_{\vv}$, let $W_{\se}\subset\Gamma$ denote the union of the face support sets of the basis functions associated with the fine grids on $\se$. Define
$$  W^{half}_{\vv}=\big(\bigcup_{\se\in{\mathcal E}_{\vv}}W_{\se}\big)\bigcap\Gamma^{half}_{\vv}.$$
Namely, $W^{half}_{\vv}$ is the intersection of $\Gamma^{half}_{\vv}$ with the union of the face fine elements adjoining $\se\in {\mathcal E}_{\vv}$.
Although the set $W^{half}_{\vv}$ looks like the {\it wire-basket} set in the BPS substructuring method, they have some differences:
$W^{half}_{\vv}$ is a vertex-related set, but the {\it wire-basket} set is subdomain-related;  the {\it wire-basket} set has zero measure in
$\Gamma$, but the set $W^{half}_{\vv}$ does not.

\begin{figure}[ht]
\centering
\begin{tabular}{ccc} 
\includegraphics[scale = 0.2]{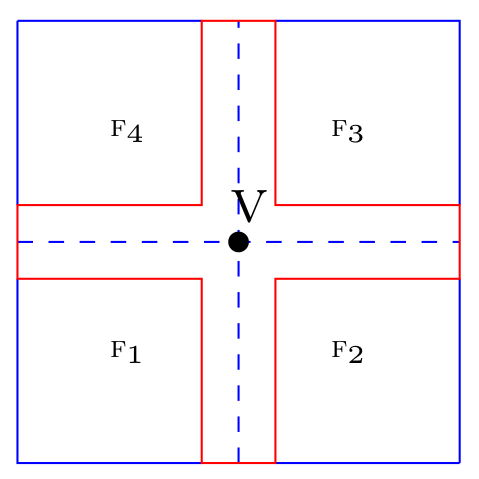}
 &\includegraphics[scale = 0.2]{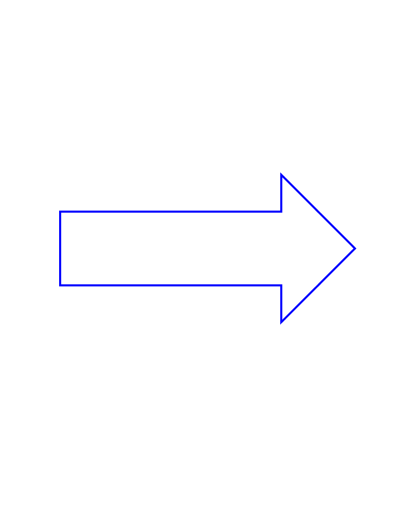}
 & \includegraphics[scale = 0.2]{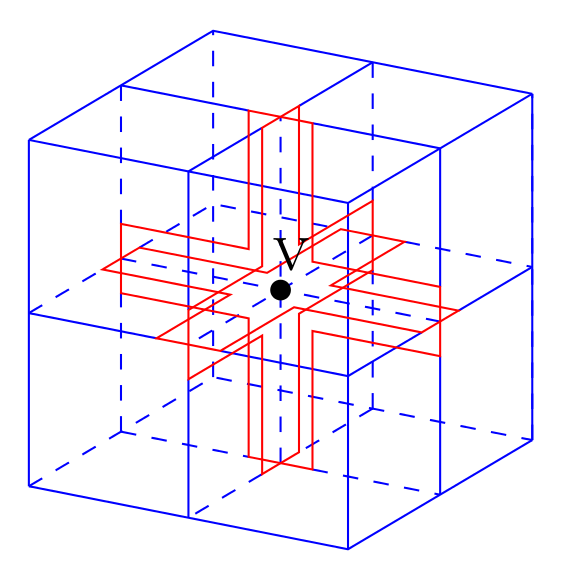}
\\
(a) 2D cross section & & (b) 3D $W^{half}_{\vv}$ and $\Omega^{half}_{\vv}$
\\
\end{tabular}
\caption{the structure of $W^{half}_{\vv}$ (red) in $\Omega^{half}_{\vv}$ (blue)}
\label{fig:crossstuc}
\end{figure}

Set
$$ F^{in}_{\vv}=\Gamma^{half}_{\vv}\backslash W^{half}_{\vv}, $$
in other words, $F^{in}_{\vv}$ is the intersection of $\Gamma^{half}_{\vv}$ with the union of the fine elements on the interior of the faces in ${\mathcal F}_{\vv}$.
Define
$$ \hat{V}^0_h(\Gamma_{\vv}^{half})=\{\v\in V^0_h(\Gamma_{\vv}^{half}):~\mbox{the~degrees~of~freedom~of}~\v~\mbox{vanish}~~\mbox{on}~F^{in}_{\vv}\}.$$
Namely, the space $\hat{V}^0_h(\Gamma_{\vv}^{half})$ keeps the ``edge" degrees of freedom of $V^0_h(\Gamma_{\vv}^{half})$ but drops the ``face" degrees of freedom.
For each $\ff\in{\mathcal F}_d$, define the local interface space
$$ V^0_h(\ff)=\{\phi\in V_h(\Gamma):~~supp~\phi\subset\ff\}.$$
It is clear that
\begin{equation}
V_h^0(\Gamma_{\vv}^{half})=\hat{V}^0_h(\Gamma^{half}_{\vv})+\big(\su_{\ff\in{\mathcal F}_{\vv}}V^0_h(\ff)\big)\cap V^0_h(\Gamma^{half}_{\vv}). \label{4.new00}
\end{equation}
Thus, by (\ref{3.new00}) we have the space decomposition
\begin{equation}
V_h(\Gamma)=\su_{\vv\in{\mathcal N}_d}\hat{V}^0_h(\Gamma^{half}_{\vv})+\su_{\ff\in{\mathcal F}_d}V^0_h(\ff). \label{3.new01}
\end{equation}

Let $\Omega_{\ff}$ be the union of $\ff$ itself and the two subdomains that have $\ff$ as their common face. Define
$$ \hat{V}^{\bot}_h(\Omega_{\vv})=\{\v\in V_h^{\bot}(\Omega_{\vv}):~~\mbox{the~trace~of}~\v~\mbox{belongs~to}~\hat{V}^0_h(\Gamma^{half}_{\vv})\}~~(\vv\in {\mathcal N}_d) $$
and
$$ V^{\bot}_h(\Omega_{\ff})=\{\v\in V_h^{\bot}(\Omega):~~\mbox{the~trace~of}~\v~\mbox{belongs~to}~V^0_h(\ff)\}~~~(\ff\in{\mathcal F}_d). $$
Namely, $\hat{V}^{\bot}_h(\Omega_{\vv})$ and $V^{\bot}_h(\Omega_{\ff})$ consist of the $A$-harmonic extensions of the functions in $\hat{V}^0_h(\Gamma^{half}_{\vv})$ and $V^0_h(\ff)$, respectively.
Corresponding to (\ref{3.new01}), we have
$$ V^{\bot}_h(\Omega)=\su_{\vv\in{\mathcal N}_d}\hat{V}^{\bot}_h(\Omega_{\vv})+\su_{\ff\in{\mathcal F}_d}V^{\bot}_h(\Omega_{\ff}). $$
Thus, by (\ref{BasicDecom}), we obtain another space decomposition
\begin{equation}
V_h(\o)=V_d(\o)+\su_{k=1}^NV_h^0(\o_k)+\su_{\vv\in {\mathcal N}_d}\hat{V}^{\bot}_h(\o_{\vv})+\su_{\ff\in{\mathcal F}_d}V_h^{\bot}(\o_{\ff}).
\label{eq:4.1}
\end{equation}

The space decomposition (\ref{eq:4.1}) seems more complicated than the space decomposition (\ref{eq:3.1}), but each subspace in the second sum and the third sum of (\ref{eq:4.1})
has different structure from $V^{\bot}_h(\o_{\vv})$, which can help us to construct cheaper local interface solvers in the next subsection.

\subsection{Approximate local interface solvers}
For $\ff\in {\mathcal F}_d$, define the operator $A_{\ff}$ as the restriction of $A$ on the subspace $V_h^{\bot}(\o_{\ff})$, i.e., it satisfies
$$ (A_{\ff}\v,\w)={\mathcal A}(\v,\w),~~\v\in V_h^{\bot}(\o_{\ff}),~\forall\w\in V_h^{\bot}(\o_{\ff}).$$
As usual, the action of $A_{\ff}^{-1}$ can be implemented by solving the following residual equation: to find $\u_{\ff}\in V_h^0(\o_{\ff})$ such that
\begin{equation}
{\mathcal A}(\u_{\ff},\w)=({\bf g},\w)-\sum_{k\in\Lambda_{\ff}}{\mathcal A}(\u_k,\w),~~~~\forall\w\in V_h^0(\o_{\ff}), \label{4.new2}
\end{equation}
where ${\bf g}\in V_h(\Omega)$ is given, and $\u_k\in V_h^0(\Omega_k)$ has been gotten by solving the local equation
$$ {\mathcal A}(\u_k,\w)=({\bf g},\w),~~~\forall\w\in V_h^0(\Omega_k).$$
But, the calculation of $\u_{\ff}$ is expensive, so we propose a new way to compute a rough approximation of $\u_{\ff}$ in the following.

Let $\Omega_{1}^{\ff}$ and $\Omega_{2}^{\ff}$ denote the subdomains sharing $\ff$ as their common face, and let $\ms{A}_{\ff}$ and $\ms{A}_{ll}$ be the stiffness
matrices generated by the basis functions on $\ff$ and in $\Omega_{l}^{\ff}$ ($l=1,2$), respectively. Then the equation (\ref{4.new2}) can be transformed into the algebraic system
\begin{equation}
\begin{pmatrix}
\ms{A}_{11} & \0& \ms{A}_{1\ff}\\
\0 & \ms{A}_{22}&\ms{A}_{1\ff}\\
\ms{A}^t_{1\ff}&\ms{A}^t_{2\ff}&\ms{A}_{\ff}
\end{pmatrix}
\begin{pmatrix}
\chi^{\ff}_{1} \\
\chi^{\ff}_{2}\\
\chi_{\ff}
\end{pmatrix}
=
\begin{pmatrix}
\0 \\
\0\\
b_{\ff}
\end{pmatrix},
\label{4.new3}
\end{equation}
where $\chi_{\ff}$ denote the dofs (i.e., coordinate vector) of $\u_{\ff}$ on $\ff$, and $b_{\ff}$ is defined by
$$ b_{\ff}=\xi_{\ff}-\ms{A}^t_{1\ff}\xi^{\ff}_1-\ms{A}^t_{2\ff}\xi^{\ff}_2, $$
with $\xi_{\ff}$ being the dofs of ${\bf g}$ on $\ff$, and $\xi^{\ff}_l$ being the dofs of $\u_l$ in $\Omega_l^{\ff}$ ($l=1,2$).

It is easy to see that the system (\ref{4.new3}) is the same as the algebraic system of the original equation (\ref{eq:2.2}) restricted in $V_h^0(\Omega_{\ff})$,
with different right hand only. Notice that, as in Step 4 of {\bf Algorithm 3.1}, the dofs. in the interior of the subdomain $\Omega^{\ff}_l$ can be gotten by computing the local harmonic extension in $\Omega^{\ff}_l$
($l=1,2$). Thus we only hope to get a rough approximation of $\chi_{\ff}$ but do not care for the accuracy of $\chi^{\ff}_{l}$ ($l=1,2$). Intuitively, the accuracy of an approximation for $\chi_{\ff}$
mainly depends on the grids nearing $\ff$ and is not sensitive to the grids far from $\ff$. Based on this observation, we can construct an auxiliary non-uniform partition $\tilde{\mathcal T}^{\ff}_{\tilde{h}}$
on $\Omega_{\ff}$, for which the original fine grids on $\ff$ are kept and the grids in $\Omega^{\ff}_{l}$ ($l=1,2$) gradually becomes coarser when nodes are far from $\ff$. Then we solve
the following auxiliary algebraic system
\begin{equation}
\begin{pmatrix}
\ms{\tilde{A}}^{\ff}_{11} & \0& \ms{\tilde{A}}_{1\ff}\\
\0 & \ms{\tilde{A}}^{\ff}_{22}&\ms{\tilde{A}}_{1\ff}\\
\ms{\tilde{A}}^t_{1\ff}&\ms{\tilde{A}}^t_{2\ff}&\ms{A}_{\ff}
\end{pmatrix}
\begin{pmatrix}
\tilde{\chi}^{\ff}_{1} \\
\tilde{\chi}^{\ff}_{2}\\
\tilde{\chi}_{\ff}
\end{pmatrix}
=
\begin{pmatrix}
\0 \\
\0\\
b_{\ff}
\end{pmatrix},
\label{4.new4}
\end{equation}
where $\ms{\tilde{A}}^{\ff}_{ll}$ denotes the stiffness matrix generated by the basis functions associated with the auxiliary grids in $\Omega^{\ff}_{l}$ ($l=1,2$). The solution
$\tilde{\chi}_{\ff}$ of the above system can be regarded as a rough
approximation of $\chi_{\ff}$. The auxiliary partition $\tilde{\mathcal T}^{\ff}_{\tilde{h}}$ (see Fig. 3) can be easily generated by the existing software \cite{gmsh}, such that the
number of the unknowns in (\ref{4.new4}) is much smaller than that in
(\ref{4.new3}), so the system (\ref{4.new4}) is very cheap to solve.


Associated with each $\hat{V}^{\bot}_h(\o_{\vv})$, we can similarly define an auxiliary partition $\tilde{\mathcal T}^{\vv}_{\tilde{h}}$ (see Fig. 3), and build the corresponding algebraic system
\begin{equation}
\begin{pmatrix}
\ms{\tilde{A}}^{\vv}_{11} & \0&\cdots&\0& \ms{\tilde{A}}_{1\vv}\\
\0 & \ms{\tilde{A}}^{\vv}_{22}&\cdots&\0&\ms{\tilde{A}}_{2\vv}\\
\vdots&\vdots&\vdots&\vdots&\vdots&\\
\0 &\cdots&\0 &\ms{\tilde{A}}^{\vv}_{mm}&\ms{\tilde{A}}^t_{m\vv}\\
\ms{\tilde{A}}^t_{1\vv}&\ms{\tilde{A}}^t_{2\vv}&\cdots&\ms{\tilde{A}}^t_{m\vv}&\ms{A}_{\vv}
\end{pmatrix}
\begin{pmatrix}
\tilde{\chi}^{\vv}_{1} \\
\tilde{\chi}^{\vv}_{2}\\
\vdots\\
\tilde{\chi}^{\vv}_{m}\\
\tilde{\chi}_{\vv}
\end{pmatrix}
=
\begin{pmatrix}
\0 \\
\0\\
\vdots\\
\0\\
b_{\vv}
\end{pmatrix},
\label{4.new5}
\end{equation}
where $\ms{A}_{\vv}$ is the stiffness matrix generated by the basis functions associated with the grids on $W^{half}_{\vv}$, and $\tilde{\chi}_{\vv}$ denotes
an approximation of the dofs. of $\u_{\vv}$ on $W^{half}_{\vv}$. Here $\u_{\vv}\in \hat{V}_h^0(\Omega_{\vv})$ is defined by
$$ {\mathcal A}(\u_{\vv},\w)=({\bf g},\w)-\sum_{k\in\Lambda_{\vv}}{\mathcal A}(\u_k,\w),~~~~\forall\w\in \hat{V}_h^0(\o_{\vv}) $$
with
$$ \hat{V}^0_h(\Omega_{\vv})=\{\v\in V_h^0(\Omega_{\vv}):~~\mbox{the~trace~of}~\v~\mbox{belongs~to}~\hat{V}^0_h(\Gamma^{half}_{\vv})\}. $$

\begin{figure}[ht]
\centering
\begin{tabular}{cc} 
\includegraphics[width =6cm]{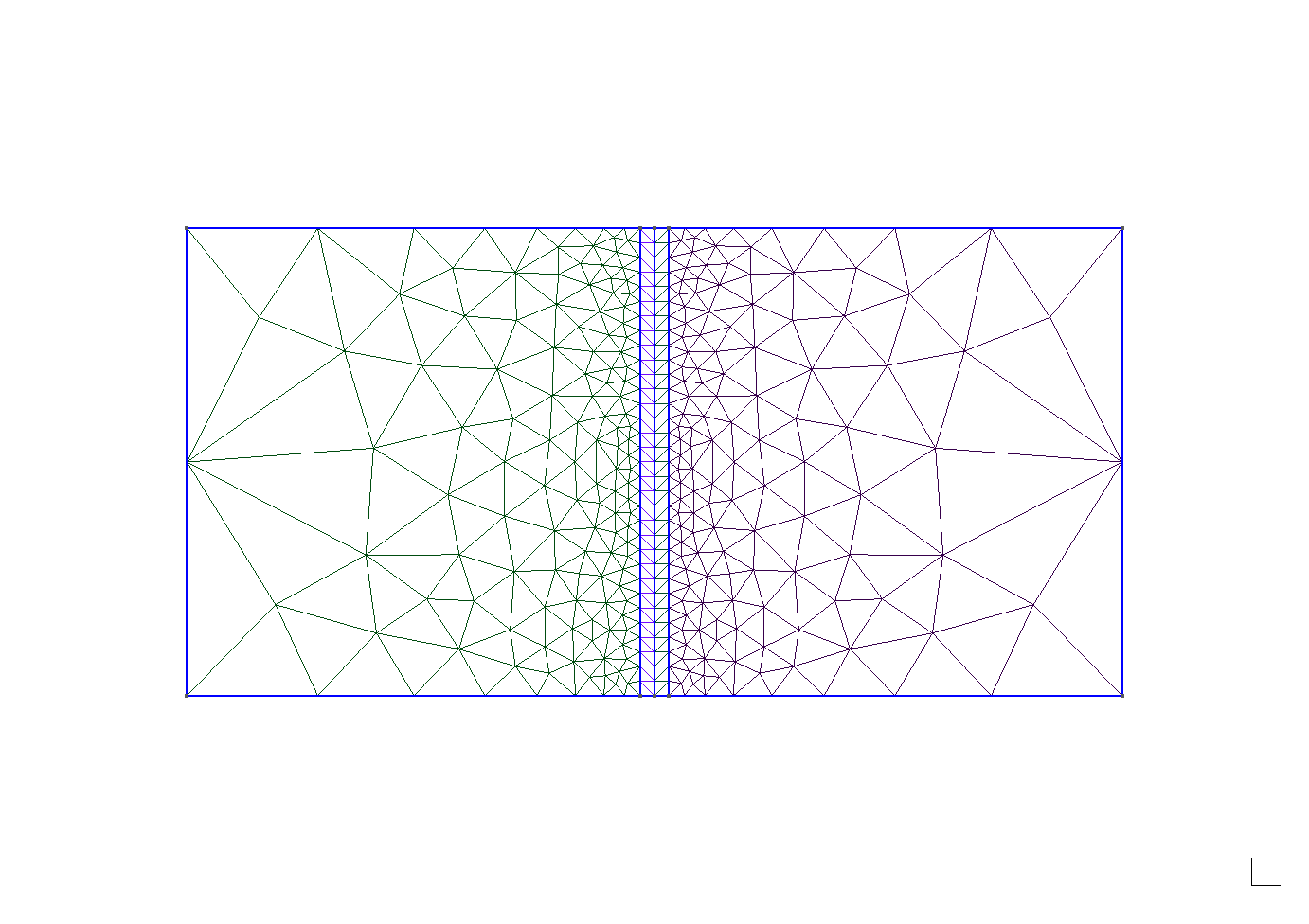}
 & \includegraphics[width = 6cm]{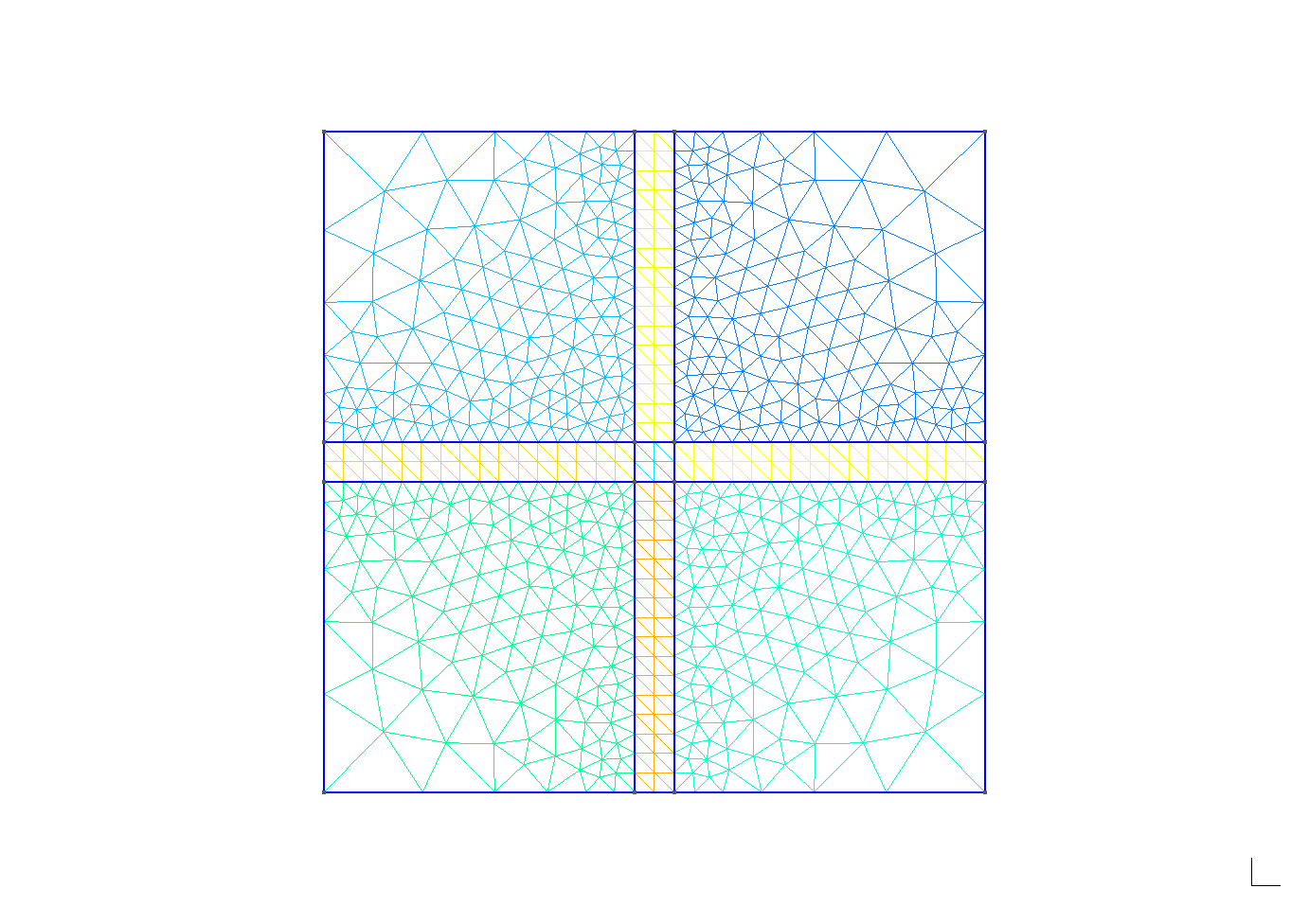}
\\
(a) $\tilde{\mathcal T}^{\ff}_{\tilde{h}}$ & (b) $\tilde{\mathcal T}^{\vv}_{\tilde{h}}$
\\
\end{tabular}
\caption{2D cross-sections of coarsening grids for $\Omega_{\ff}$ and $\Omega^{half}_{\vv}$}
\label{fig:coarse_2dmesh}
\end{figure}


According to the above discussions, the approximate solvers $\tilde{A}_{\ff}: V_h^{\bot}(\o_{\ff})\rightarrow V_h^{\bot}(\o_{\ff})$ and $\tilde{A}_{\vv}: \hat{V}_h^{\bot}(\o_{\vv})\rightarrow \hat{V}_h^{\bot}(\o_{\vv})$
can be defined as follows: for ${\bf g}\in V_h^{\bot}(\o_{\ff})$, we define
$\tilde{\u}_{\ff}=\tilde{A}^{-1}_{\ff}{\bf g} \in V_h^{\bot}(\o_{\ff})$ and $\tilde{\u}_{\vv}=\tilde{A}^{-1}_{\vv}{\bf g} \in \hat{V}_h^{\bot}(\o_{\vv})$, such that the dofs of $\tilde{\u}_{\ff}$ on $\ff$
and $W^{half}_{\vv}$ equal to $\tilde{\chi}_{\ff}$ and $\tilde{\chi}_{\vv}$ computed by solving
(\ref{4.new4}) and (\ref{4.new5}), respectively.

\begin{remark} The construction of the above approximate solvers do not depend on the original bilinear form ${\mathcal A}(\cdot,\cdot)$ and the interior grids of the subdomains for coarsening.
Thus we need not to estimate the norm induced from ${\mathcal A}(\cdot,\cdot)$ on $\ff$ and $W_{\vv}^{half}$, and need not to give special assumptions on the original grids.
\end{remark}

\begin{remark} For the coarsening technique introduced in this subsection, the grids on the considered local interface with non-zero dofs need to be kept, so the number of nodes on such an interface should be much smaller
than that on the boundary of every subdomain for coarsening, otherwise, the number of nodes of the coarsening partition is still great (see the data in Subsection 6.4). Because of this, we did not make the coarsening directly
for the original space $V^{\bot}_h(\o_{\vv})$ defined in Subsection 3.1, and we have to build the decomposition (\ref{4.new00}) and define the spaces $\hat{V}^{\bot}_h(\o_{\vv})$ and $V^{\bot}_h(\o_{\ff})$. Of course,
if there is no degree of freedom on the edges in ${\mathcal E}_{\vv}$ (Raviart-Thomas elements), then we need not to define the space $\hat{V}^{\bot}_h(\o_{\vv})$ (which only grasps the degrees of freedom on coarse edges). For this case, the second sum in (\ref{eq:4.1}) will not appear.
\end{remark}


\subsection{Preconditioner}
By using the space decomposition in Subsection 4.1 and the approximate interface solvers in Subsection 4.2, we can define the second preconditioner for $A$ as
\be
B_{II}^{-1}=B_d^{-1}Q_d+\su_{k=1}^NA^{-1}_k Q_k+\su_{\vv\in {\mathcal
N}_d}\tilde{A}_{\vv}^{-1}Q_{\vv}+
\sum_{\ff\in{\mathcal F}_d} \tilde{A}_{\ff}^{-1}Q_{\ff}.\l{precon2} \ee
When there is no degree of freedom on the edges in ${\mathcal E}_{\vv}$ (Raviart-Thomas elemens),
the second sum in the preconditioner would be dropped.

The action of $B_{II}^{-1}$ can be described by the following algorithm.\\

{\bf Algorithm 4.1}. For ${\bf g}\in V_h(\o)$, we can compute
$\u=B_{II}^{-1}{\bf g}$ in five steps.

Step 1. Solve the system of $\u_d\in V_d(\o)$:
$$ {\mathcal A}(\u_d, \v_d)=({\bf g},\v_d),~~~\forall\v_d\in V_d(\o); $$

Step 2. Solve the following systems of $\u_k\in V_h^0(\o_k)$ in parallel:
$$ {\mathcal A}(\u_k, \v)=({\bf g},\v),~~~\forall\v\in V^0_h(\o_k)~~~(k=1,\cdots,N);
$$

Step 3. Solve the system (\ref{4.new4}) to get the dofs $\tilde{\chi}_{\ff}$ of $\tilde{\u}_{\ff}$ on $\ff$ in parallel for every $\ff\in{\mathcal F}_{\vv}$;

Step 4. Solve the system (\ref{4.new5}) to get the dofs $\tilde{\chi}_{\vv}$ of $\tilde{\u}_{\vv}$ on $W^{half}_{\vv}$ in parallel for every $\vv\in{\mathcal N}_{d}$;

Step 5. Use all $\tilde{\chi}_{\ff}$ and $\tilde{\chi}_{\vv}$ to get the dofs. $\tilde{\chi}_{\partial\Omega_k}$ of $\su_{\ff\subset\partial\Omega_k}\tilde{\u}_{\ff}+\su_{\vv\in \partial\Omega_k}\tilde{\u}_{\vv}$ on $\partial\Omega_k$. Compute the discrete $A$-extension $\u^{\bot}_k$ in parallel, such that $\u^{\bot}_k$ has the dofs $\tilde{\chi}_{\partial\Omega_k}$ on $\partial\Omega_k$ and satisfies
$$ {\mathcal A}(\u^{\bot}_k,\v)=0,~~~~\forall\v\in V_h^0(\Omega_k)~~~(k=1,\cdots,N).$$
Finally, we define
$$ \u=\u_d+\su_{k=1}^N(\u_k+\u^{\bot}_k). $$

\begin{remark} The essential difference between the proposed substructuring method and the existing substructuring methods is that novel local interface solvers are used in Step 3 and Step 4 of the above algorithm.
As explained in Subsection 4.2, the interface solvers are cheap and easy to implement.
In fact, the auxiliary subproblems needed to be solved in Step 3 and Step 4 have very small dofs. (see the data in Section 6). When there is no degree of freedom on the edges in ${\mathcal E}_{\vv}$ (Raviart-Thomas elemens),
we need not to implement Step 4. Notice that the local problems in Step 5 have the same stiffness matrices with that in Step 2
(with different right hands only). Thus the implementation of Step 5 is very cheap by using LU decomposition made in Step 2 for each local stiffness matrix.
\end{remark}

\section{Applications}
\setcounter{equation}{0}
 In this section, we
introduce two typical elliptic-type equations.

\subsection{Linear elasticity problems}
Let's consider the linear elasticity problem:
\be \begin{cases}
-\sum\limits_{j=1}^3\frac{\partial \sigma_{ij}}{\partial x_{j}}(\bm{u})=f_{i},\q in~~\O \\
 \q\q\q\q\q \bm{u}=0, \q on~~\partial\Omega \end{cases}\l{eq:5.1}
\ee
where $\bm{f}=(f_1~f_2~f_3)^{T}$ is an internal volume force, e.g. gravity (cf. \cite{ChenH2011}). The linearized strain tensor is defined by
\[\varepsilon = \varepsilon(\bm{u}) = [\varepsilon_{ij} = \frac{1}{2}(\frac{\partial
u_{i}}{\partial x_{i}} +\frac{\partial u_{j}}{\partial x_{i}})] \]
and
\[ \sigma_{ij}(\bm{u}) := \lambda\delta_{ij}div\bm{u} + 2\mu\varepsilon_{ij}, \]
where $\lambda$ and $\mu$ are the $Lam\acute{e}$ parameters (cf. \cite{re25}), which are positive
functions.

As usual, let $H_{0}^1(\Omega)\subset H^1(\Omega)$ denote the space consisting of the
functions having the zero trace on $\partial \Omega$. We introduce the vector
value Sobolev space $(H_{0}^1(\Omega))^3$, equipped with the usual product norm as follows:
\[||\bm{u}||_{1,\Omega}:=(|\bm{u}|_{H^1(\Omega)}^2 + ||\bm{u}||_{L_{2}(\Omega)}^2)^{\frac{1}{2} }\]
with $||\bm{u}||_{L_{2}(\Omega)}^2:=\int_{\Omega}|\bm{u}|^2dx$ and $|\bm{u}|_{H^1(\Omega)}^2:=||\nabla \bm{u}||_{L_{2}(\Omega)}^2$.
Concerning the
variational problem (\ref{eq:2.1}), we have $V(\Omega):=[H_{0}^1(\Omega)]^3$,
\[ \bm{\A}(\bm{u},\bm{v}) = \int_{\Omega}(2\mu\varepsilon(\bm{u}):\varepsilon(\bm{v}) + \lambda div\bm{u}\cdot div\bm{v})dx \]
and
\[ \langle \bm{F},\bm{v} \rangle = \int_{\Omega}\bm{f}\cdot\bm{v}dx \]
with
\[\varepsilon(\bm{u}):\varepsilon(\bm{v}) := \sum_{i,j=1}^n \varepsilon_{ij}(\bm{u})\varepsilon_{ij}(\bm{v}).\]

Let $R(K)$ be a subset of all linear polynomials on the element
$K$ of the form:
$$ R(K)=\Big\{\bm{A}\cdot\x+\bm{C};~ \bm{A}\in \mathds{R}^{3\times3}, \bm{C} \in \mathds{R}^3, ~\x\in K\Big\} . $$

Assume that $\O$ can be written as the union of polyhedral subdomains $D_1$, $\cdots$, $D_{N_0}$, such that $\lambda(x) =\lambda_r$ and $\mu(x) = \mu_r$
for $x\in D_r$, with $\lambda_r$ and $\mu_r$ being positive constants. In applications, $N_0$ is a {\it fixed} positive integer, so the diameter of each $D_r$ is $O(1)$.
It is certain that the subdomains $\Omega_k$ should satisfy the condition: each $D_r$ is the union of some subdomains in $\{\Omega_k\}$.

\subsection{Maxwell's equations}
For the time-dependent Maxwell's equations, we need to solve the following {\bf curlcurl}-system at each time step (see \cite{Cessenat1998, Hiptmair2002, Monk2003}):
\begin{equation}
\left \{
\begin{array}{rrr}
\c(\al\,\c\,\u)+\beta\u=\f, &in& \q\o,
\\
\u \times \bm{n} = 0, &on& \partial\Omega
\end{array}
\right.
\label{eq:5.2}
\end{equation}
where the coefficients $\al(\x)$ and $\beta(\x)$ are two positive
bounded functions in $\o$, and $\n$ is the unit outward normal vector on $\p\o$.

Let $H({\bf curl};\o)$ be the Sobolev space consisting of all square
integrable functions whose {\bf curl}'s are also square integrable
in $\o$, and $ H_0({\bf curl};\o)$ be the subspace of $H({\bf
curl};\o)$ of all functions whose tangential components vanishing
on $\p\o$. In order to get the weak form of
(\ref{eq:5.2}), just like linear elasticity problems, we define $V(\Omega)=H_0({\bf curl})$,
\[ \bm{\A}(\bm{u},\bm{v}) = \int_{\Omega}(\alpha~\c~\u \cdot \c~\v
+ \beta~\u \cdot \v)dx \]
and
\[ \langle \bm{F},\bm{v} \rangle = \int_{\Omega}\f\cdot\bm{v}dx. \]

Let $R(K)$ be a subset of all linear polynomials on the element
$K$ of the form:
$$ R(K)=\Big\{\a+\b\ti\x;~ \a,\b\in \mathds{R}^3, ~\x\in K\Big\} \,. $$
It is well-known that for any $\v\in V_{h}(\o)$, its tangential
components are continuous on all edges of each element in the
triangulation ${\mathcal T}_h$. Moreover, each edge element function
$\v$ in $V_{h}(\o)$ is uniquely determined by its moments on each edge
$e$ of ${\mathcal T}_h$:
\begin{equation}
 \Big\{\la_e(\v)=\int_e\v\cdot\t_e ds; ~ e\in {\mathcal E}_h\Big\},
 \label{eq:5.3}
\end{equation}
where $\t_e$ denotes the unit vector on the edge $e$.

As in the last subsection, we assume that $\O$ can be written as the union of polyhedral
subdomains $D_1$, $\cdots$, $D_{N_0}$ with $N_0$ being a fixed positive integer, such that
$\alpha(x)= \alpha_r$ and $\beta(x)=\beta_r$ for $x\in D_r$, where
every $\alpha_r$ and $\beta_r$ is a positive constant. Let the subdomains $\Omega_k$ satisfy the condition: each $D_r$ is the union of some subdomains in $\{\Omega_k\}$.

\section{Numerical Experiments}
\setcounter{equation}{0} In this section, we report some numerical results to illustrate the effectiveness of the proposed substructuring
preconditioners.

We consider the models introduced in Section 5, with $\Omega=(0,1)^3$, and we make tests for different distributions of the coefficients $\lambda(x)$, $\mu(x)$, $\alpha(\x)$ and $\beta(\x)$:

\medskip
\noindent{\bf Case (i)}: the coefficients have no jump, i.e., $\lambda(x) = \mu(x) = 1$ (linear elasticity problems) or $\alpha(\x)=\beta(\x)=1$
(Maxwell's equations).

\medskip
\noindent{\bf Case (ii)}: the coefficients have large jumps, i.e.,
\[\lambda(\x)=\left\{
\begin{array}{l}
\lambda_0,\quad\quad in~D\cr  \\
~~1,~\quad in~\o\backslash D,
\end{array}
\right. \quad\quad  \mu(\x)=\left\{
\begin{array}{l}
\mu_0,\quad\quad in~D,\cr  \\
~~1,~\quad in~\o\backslash D
\end{array}
\right.
\]
for linear elasticity problems and
\[\alpha(\x)=\left\{
\begin{array}{l}
\al_0,\quad\quad in~D\cr  \\
~~1,~\quad in~\o\backslash D,
\end{array}
\right. \quad\quad  \beta(\x)=\left\{
\begin{array}{l}
\beta_0,\quad\quad in~D,\cr  \\
~~1,~\quad in~\o\backslash D
\end{array}
\right.
\]
for Maxwell's equations. Here $D\subset\o$ is a union of several subdomains $\o_k$. We consider two
choices of $D$:

$$ \mbox{Choice (1)}. ~~~D=[\dfrac{1}{4},~\frac{1}{2}]^3;~~~~\mbox{Choice (2)}.~~~ D=[\frac{1}{4},~\frac{1}{2}]^3\bigcup[\frac{1}{2},~\frac{3}{4}]^3. $$

In our experiments, we define domain decomposition and finite element partition as follows.
At first, we divide the domain into $n^3$ smaller cubes $\Omega_1$, $\Omega_2\cdots\Omega_N$, which have the same length $d$ of edges, i.e., $d = 1/n$.
We require that $D\subset\o$ is just the union of some subdomains in $\{\o_k\}$, which yields the desired domain decomposition. Next, we divide each subdomain $\Omega_k$ into $m^3$
fine cubes, with the same size $h = 1/(mn)$. All the fine cubes constitute a partition ${\mathcal T}_h$ consisting of hexahedral elements.
If we further divide each fine cube into 5 or 6 tetrahedrons in the standard way, then all the generated tetrahedrons constitute a partition ${\mathcal T}_h$ consisting of
tetrahedral elements.

We discretize the models by the linear finite element methods, and we apply the PCG method with the proposed preconditioners to solve the resulting algebraic systems. The PCG iteration is terminated
in our experiments when the relative residual is less than $10^{-6}$. We will report the iteration counts in the rest of this section.

\subsection{Tests for linear elasticity problems}
In this subsection, we consider an example of the linear elasticity
problem. We choose the right-hand side $\bm{f}$ of system (\ref{eq:5.1}) such that the analytic solution $\u=(u_1, u_2, u_3)^T$ is given by:
\begin{eqnarray*}
 u_1 &=& x(x-1)y(y-1)z(z-1)\\
 u_2 &=& x(x-1)y(y-1)z(z-1)\\
 u_3 &=& x(x-1)y(y-1)z(z-1)
 \end{eqnarray*}
where the coefficients $\lambda(x) = \mu(x) = 1$. In our experiments, the
right-hand side $\bm{f}$ is fixed.

\subsubsection{Efficiency of the first preconditioner}
In this part, we test the action of the preconditioner $B_I$ described by {\bf
Algorithm 3.1}. We use both hexahedral partition and tetrahedral partition in
our experiments. We first consider the case of hexahedral partition.
The iteration counts of the PCG method with $B_I$ are listed in  Table \ref{E-hp-I-sm} (for {\bf Case} (i)) and  Table \ref{E-hp-I-nsm} (for {\bf Case} (ii)).

\vskip 0.1in
  \begin{center}  
 \tabcaption{}\label{E-hp-I-sm}
   Iteration counts of PCG with the preconditioner $B_I$ (hexahedral elements): the coefficients have no jumps
  \vskip 0.2in
\begin{tabular}{|c|c|c|c|c|}
\hline
$ m \backslash n$ &  ~~4~~ &  ~~6~~ &  ~~8~~ & ~~10~~ \\\hline
                 4&  15& 15& 15& 14 \\\hline
                 8&  16& 16& 16& 16 \\\hline
                16&  18& 18& 19& 19 \\\hline
                32&  21& 21& 21& 21   \\\hline
\end{tabular}
\end{center}
 \vskip 0.1in

\vskip 0.1in
\begin{center}
 \tabcaption{}\label{E-hp-I-nsm}
  Iteration counts of PCG with the preconditioner $B_I$ (hexahedral elements): the coefficients have large jumps
  \vskip 0.2in
\begin{tabular}{|c|c|c|c|c|c|c|c|c|}
\hline
$\mbox{}$
&\multicolumn{4}{c|}{Choice~(1)~of~$D$}
&\multicolumn{4}{c|}{Choice~(2)~of~$D$}\\\hline
$\mbox{}$
&\multicolumn{2}{c|}{$ \lambda_0 =\mu_0=10^{-5}$}
&\multicolumn{2}{c|}{$ \lambda_0 =\mu_0=10^{5} $}
&\multicolumn{2}{c|}{$ \lambda_0 =\mu_0=10^{-5} $}
&\multicolumn{2}{c|}{$ \lambda_0 =\mu_0=10^{5}  $}\\\hline
$~m ~ \backslash ~ n$& ~~4~~&~~8~~&~~4~~&~~8~~&~~4~~&~~8~~&
~~4~~&~~8~~\\\hline
 8 &14    &16  &19  &19 &14 &16 &18  &19  \\\hline
 16&16    &19  &22  &21 &16 &19 &21  &21   \\\hline
 24&18    &21  &23  &23 &18 &21 &23  &23  \\\hline
 32&19    &22  &26  &24 &18 &22 &24  &24     \\\hline
\end{tabular}
\end{center}
 \vskip 0.1in

We observe from Table \ref{E-hp-I-sm} that, when the coefficients is smooth, the iteration counts of PCG method grows slowly when $m =d/h$ increases but $n=1/d$ is fixed,
and almost unchange when $m$ is fixed but $n$ increases. The data in Table \ref{E-hp-I-nsm} indicate that, even if the coefficients have large jumps, the iteration counts of PCG still grows slowly.
It confirms that the preconditioner $B_I$ is effective for the system
 arising from nodal element discretization for linear elasticity problems.

Next we consider the case with tetrahedral partition. Since the subdomains are hexahedrons, the coarse space associated with the subdomains is not a subspace of
the fine tetrahedral element space. Because of this, we further divide each cubic subdomain into 5 or 6 tetrahedrons, and use all the tetrahedral subdomains to
define a nested coarse space. Notice that the resulting tetrahedral coarse space has the same number of the degrees of freedom as the original hexahedral coarse
space, i.e., this change will not increase the cost for implementing the coarse solver.

We list the iteration counts of the PCG method with $B_I$ in Table \ref{E-tp-I-sm} (for
{\bf Case} (i)) and Table \ref{E-tp-I-nsm} (for {\bf Case} (ii)).

\vskip 0.1in
  \begin{center}  
  \tabcaption{}\label{E-tp-I-sm}
  Iteration counts of PCG with the preconditioner $B_I$ (tetrahedral elements): the coefficients have no jump
  \vskip 0.2in
\begin{tabular}{|c|c|c|c|c|}
\hline
$ m \backslash n$  &  ~~4~~ &  ~~6~~ &  ~~8~~ & ~~10~~ \\\hline
                 8&   20& 20& 20& 19 \\\hline
                16&   23& 23& 22& 21 \\\hline
                24&   24& 24& 23& 23  \\\hline
                32&   25& 25& 24& 24   \\\hline
\end{tabular}
\end{center}
 \vskip 0.1in

\vskip 0.1in
\begin{center}
  \tabcaption{}\label{E-tp-I-nsm}
   Iteration counts of PCG with the preconditioner $B_I$ (tetrahedral elements): the coefficients have large jumps
  \vskip 0.2in
\begin{tabular}{|c|c|c|c|c|c|c|c|c|}
\hline
$\mbox{}$
&\multicolumn{4}{c|}{Choice~(1)~of~$D$}
&\multicolumn{4}{c|}{Choice~(2)~of~$D$}\\\hline
$\mbox{}$
&\multicolumn{2}{c|}{$ \lambda_0 =\mu_0=10^{-5}$}
&\multicolumn{2}{c|}{$ \lambda_0 =\mu_0=10^{5} $}
&\multicolumn{2}{c|}{$ \lambda_0 =\mu_0=10^{-5} $}
&\multicolumn{2}{c|}{$ \lambda_0 =\mu_0=10^{5}  $}\\\hline
$~m ~ \backslash ~ n$& ~~4~~&~~8~~&~~4~~&~~8~~&~~4~~&~~8~~&
~~4~~&~~8~~\\\hline
 8 &17    &20  &27  &23 &17 &21 &27  &23  \\\hline
 16&20    &22  &29  &25 &20 &24 &29  &25   \\\hline
 24&21    &24  &31  &27 &21 &25 &31  &27     \\\hline
 32&22    &25  &32  &28 &22 &26 &32  &28  \\\hline
\end{tabular}
\end{center}
 \vskip 0.1in

We observe that the iteration counts of PCG in these two tables vary
stably for the considered two cases (even if the coefficents have large jumps). In addition,
we can see that the convergence rate of PCG is same as in the case with the hexahedron elements.

\subsubsection{Efficiency of the second preconditioner}
 In this subsection we investigate the efficiency of the preconditioner $B_{II}$ described by {\bf Algorithm 4.1}.

Firstly, we consider the case of hexahedral elements.
The iteration counts of the PCG method with $B_{II}$ are listed in
 Table \ref{E-hp-II-sm} (for {\bf Case} (i)) and Table \ref{E-hp-II-nsm} (for {\bf Case} (ii)).

\vskip 0.1in
\begin{center}
\tabcaption{}\label{E-hp-II-sm}
  Iteration counts of PCG with the preconditioner $B_{II}$ (hexahedral elements): the coefficients have no jump
 \vskip 0.2in
\begin{tabular}{|c|c|c|c|c|}
\hline
$ m \backslash n$ &  ~~4~~  & ~~6~~  &~~8~~   &~~10~~ \\\hline
 8& 19 &19 &19  &19  \\\hline
16& 22 &22 &22  &22 \\\hline
24& 24 &24 &24  &23 \\\hline
32& 25 &25 &24  &24  \\\hline
\end{tabular}
\end{center}
 \vskip 0.1in

\vskip 0.1in
\begin{center}
\tabcaption{}\label{E-hp-II-nsm}
  Iteration counts of PCG with the preconditioner $B_{II}$ (hexahedral elements): the coefficients have large jumps
\vskip 0.2in
\begin{tabular}{|c|c|c|c|c|c|c|c|c|}\hline
$\mbox{}$
&\multicolumn{4}{c|}{Choice~(1)~of~$D$}
&\multicolumn{4}{c|}{Choice~(2)~of~$D$}\\\hline
$\mbox{}$
&\multicolumn{2}{c|}{$ \lambda_0 =\mu_0=10^{-5}$}
&\multicolumn{2}{c|}{$ \lambda_0 =\mu_0=10^{5} $}
&\multicolumn{2}{c|}{$ \lambda_0 =\mu_0=10^{-5} $}
&\multicolumn{2}{c|}{$ \lambda_0 =\mu_0=10^{5}  $}\\\hline
$~m ~ \backslash n$& ~~4~~&~~8~~&~~4~~&~~8~~&~~4~~&~~8~~&
~~4~~&~~8~~\\\hline
 8 &17    &20  &24  &23 &17 &20 &25  &23     \\\hline
 16&20    &24  &28  &27 &20 &24 &28  &27     \\\hline
 24&22    &26  &30  &29 &22 &26 &31  &29     \\\hline
 32&23    &27  &31  &30 &23 &27 &32  &30      \\\hline
\end{tabular}

\end{center}
 \vskip 0.1in

From Table \ref{E-hp-II-sm}, we observe that the rate of convergence of PCG with $B_{II}$ is same as that with $B_I$. In addition, we found that the iteration counts in Table \ref{E-hp-II-sm} are slightly more than
that in  Table \ref{E-hp-I-sm} when the values of $m, n$ are same in these two tables. But in
each PCG iteration step, the calculation of $B^{-1}_{II}g$ is much cheaper
than $B^{-1}_{I}g$ when $d/h$ is large enough (we will investigate this question in the final subsection of this section). We can see from
 Table \ref{E-hp-II-nsm} that, even if coefficients have large jumps, the iteration counts vary stably. This means that the preconditioner $B_{II}$ is not only cheaper, but also effective for elasticity
problems.

Next we consider the case of tetrahedral partition. Here we construct a
coarse space as in the last subsection for $B_I$.
The iteration counts of the PCG are listed in
Table \ref{E-tp-II-sm} (for {\bf Case} (i)) and
 Table \ref{E-tp-II-nsm} (for {\bf Case} (ii)).

\vskip 0.1in
\begin{center}
\tabcaption{}\label{E-tp-II-sm}
  Iteration counts of PCG with the preconditioner $B_{II}$ (tetrahedral elements): the coefficients have no jump
 \vskip 0.2in
\begin{tabular}{|c|c|c|c|c|}
\hline
$ m \backslash n$ &  ~~4~~  & ~~6~~  &~~8~~   &~~10~~ \\\hline
 8& 22 &22 &21  &21  \\\hline
 16&25 &25 &24  &24 \\\hline
 24&27 &27 &26  &25 \\\hline
 32&29 &28 &27  &27  \\\hline
\end{tabular}

\end{center}
 \vskip 0.1in

\vskip 0.1in
\begin{center}
\tabcaption{}\label{E-tp-II-nsm}
  Iteration counts of PCG with the preconditioner $B_{II}$ (tetrahedral elements): the coefficients have large jumps
\vskip 0.2in
\begin{tabular}{|c|c|c|c|c|c|c|c|c|}\hline
$\mbox{}$
&\multicolumn{4}{c|}{Choice~(1)~of~$D$}
&\multicolumn{4}{c|}{Choice~(2)~of~$D$}\\\hline
$\mbox{}$
&\multicolumn{2}{c|}{$ \lambda_0 =\mu_0=10^{-5}$}
&\multicolumn{2}{c|}{$ \lambda_0 =\mu_0=10^{5} $}
&\multicolumn{2}{c|}{$ \lambda_0 =\mu_0=10^{-5} $}
&\multicolumn{2}{c|}{$ \lambda_0 =\mu_0=10^{5}  $}\\\hline
$~m ~ \backslash n$& ~~4~~&~~8~~&~~4~~&~~8~~&~~4~~&~~8~~&
~~4~~&~~8~~\\\hline

 8 &19    &22  &31  &26 &19 &23 &30  &27     \\\hline
 16&23    &27  &35  &30 &22 &27 &34  &31     \\\hline
 24&24    &29  &37  &32 &24 &29 &37  &33     \\\hline
 32&26    &31  &39  &34 &26 &31 &39  &35       \\\hline
\end{tabular}
\end{center}
 \vskip 0.1in

Like the case of hexahedral elements, the preconditioner $B_{II}$ is still effective for the case of tetrahedral elements.

\subsection{Tests for Maxwell's equations}
In this subsection, we consider Maxwell's equations.
Let the right-hand side $\bm{f}$ in the equations
(\ref{eq:5.2}) to be selected such that
the exact solution $\u=(u_1, u_2, u_3)^T$ is given by
\begin{eqnarray*}
 u_1&=&xyz(x-1)(y-1)(z-1)\,, \\
 u_2&=& \sin (\pi x) \sin(\pi y) \sin(\pi z)\,, \\
 u_3&=& (1-e^x)(1-e^{x-1})(1-e^y)(1-e^{y-1})(1-e^z)(1-e^{z-1})\,,
\end{eqnarray*}
where the coefficients $\alpha(\x)$ and $\beta(\x)$ are both constant $1$.
This right-hand side $\bm{f}$ is also fixed in our experiments.

\subsubsection{Efficiency of the first preconditioner}
In this part, we investigate the effectiveness of the preconditioner $B_I$ described by {\bf
Algorithm 3.1}. We first consider the case of hexahedral elements.
The iteration counts of the PCG method with $B_I$ are listed in
 Table \ref{M-hp-I-sm} (for {\bf Case} (i)) and  Table \ref{M-hp-I-nsm} (for {\bf Case} (ii)).

\vskip 0.1in
\begin{center}
\tabcaption{}\label{M-hp-I-sm}
 Iteration counts of PCG with the preconditioner $B_{I}$ (hexahedral elements): the coefficients have no jump
\vskip 0.2in
\begin{tabular}{|c|c|c|c|c|}
\hline
$ m \backslash n$ & ~~4~~  & ~~6~~  &~~8~~   &~~10~~ \\\hline
 8& 16  & 15 & 15  &15    \\\hline
 16&17  & 18 & 18  &17   \\\hline
 24&19  & 19 & 19  &18   \\\hline
 32&20  & 20 & 20  &20     \\\hline
\end{tabular}
 \end{center}
  \vskip 0.1in

\vskip 0.1in
\begin{center}
\tabcaption{}\label{M-hp-I-nsm}
 Iteration counts of PCG with the preconditioner $B_{I}$ (hexahedral elements): the coefficients have large jumps
\vskip 0.2in
\begin{tabular}{|c|c|c|c|c|c|c|c|c|}\hline
$\mbox{}$& \multicolumn{4}{c|}{Choice~(1)~of~$D$}&
\multicolumn{4}{c|}{Choice~(2)~of~$D$}\\\hline $\mbox{}$&
\multicolumn{2}{c|}{$ \beta_0 =\alpha_0 =10^{-5}  $
}&\multicolumn{2}{c|}{$ \beta_0 =\alpha_0=10^{5} $}&
\multicolumn{2}{c|}{$ \beta_0 =\alpha_0=10^{-5}  $
}&\multicolumn{2}{c|}{$ \beta_0 =\alpha_0=10^{5} $}\\\hline
$~m ~ \backslash n$& ~~4~~&~~8~~&~~4~~&~~8~~&
~~4~~&~~8~~&~~4~~&~~8~~\\\hline

 8& 13    &15   &19  &17     &13   &15  &19  &19    \\\hline
 16&15    &17   &21  &20     &15   &17  &22  &22  \\\hline
 24&16    &18   &23  &21     &16   &19  &24  &24    \\\hline
 32&16    &19   &24  &22     &16   &19  &25  &25    \\\hline
\end{tabular}

\end{center}
 \vskip 0.1in

We observe from the above two tables that, although the coarse
space is chosen as the simplest one for Maxwell's
equations, the iteration counts of the PCG method with the preconditioner $B_I$ grow logarihmically with $m=d/h$ only, not depend on $n=1/d$, even if the coefficients have large jumps.

Now we consider the case of tetrahedral elements. For this case, we can not simply consider the coarse space corresponding to the hexahedral subdomain partition. If we
divide each hexahedral subdomain into 5 or 6 tetrahedral subdomains and use the tetrahedral coarse space as in the last section, then the tetrahedral coarse space
have much more degrees of freedom than the natural hexahedral coarse space, since the degrees of freedom are defined on the coarse edges for Maxwell's equations. A natural idea is to
define a tetrahedral coarse space as the image space of the interpolation operator acting on the natural hexahedral coarse space. Then the degrees of freedom are not increased in the resulting tetrahedral coarse space.
The iteration counts of the PCG method are listed in  Table \ref{M-tp-I-sm} (for {\bf Case} (i)) and  Table \ref{M-tp-I-nsm} ( {\bf Case} (ii)).

\vskip 0.1in
\begin{center}
\tabcaption{}\label{M-tp-I-sm}
 Iteration counts of PCG with the preconditioner $B_{I}$ (tetrahedral elements): the coefficients have no jump
\vskip 0.2in
\begin{tabular}{|c|c|c|c|c|c|}
\hline
$ m \backslash n$ & ~~4~~  & ~~6~~  &~~8~~   &~~10~~ \\\hline
  8  & 18 & 17 &17 &16 \\\hline
 16  & 20 &20  &19 &19\\\hline
 24  & 22 &21  &21 &21  \\\hline
 32  & 22 &22  &22 &21  \\\hline
\end{tabular}
 \end{center}
  \vskip 0.1in

\vskip 0.1in
\begin{center}
\tabcaption{}\label{M-tp-I-nsm}
 Iteration counts of PCG with the preconditioner $B_{I}$ (tetrahedral elements): the coefficients have large jumps
\vskip 0.2in
\begin{tabular}{|c|c|c|c|c|c|c|c|c|}\hline
$\mbox{}$& \multicolumn{4}{c|}{Choice~(1)~of~$D$}&
\multicolumn{4}{c|}{Choice~(2)~of~$D$}\\\hline $\mbox{}$&
\multicolumn{2}{c|}{$ \beta_0 =\alpha_0 =10^{-5}  $
}&\multicolumn{2}{c|}{$ \beta_0 =\alpha_0=10^{5} $}&
\multicolumn{2}{c|}{$ \beta_0 =\alpha_0=10^{-5}  $
}&\multicolumn{2}{c|}{$ \beta_0 =\alpha_0=10^{5} $}\\\hline
$~m ~ \backslash n$& ~~4~~&~~8~~&~~4~~&~~8~~&
~~4~~&~~8~~&~~4~~&~~8~~\\\hline
 8& 15    &17  &21  &20     &15    &17  &21  &21    \\\hline
 16&17    &20  &24  &22     &17    &20  &24  &24    \\\hline
 24&18    &21  &25  &23     &18    &21  &26  &26      \\\hline
 32&19    &22  &26  &24     &19    &22  &27  &28      \\\hline
\end{tabular}

\end{center}
 \vskip 0.1in

From the above two tables, we can see that the iteration counts vary
stably and the PCG iteration has the same convergence rate as in the case of the linear elasticity problem.

\subsubsection{Efficiency of the second preconditioner}
In this part, we investigate the effectiveness $B_{II}$ for the case of Maxwell's equations.

Firstly, we consider the hexahedral partition. We list the iteration counts of PCG method in  Table \ref{M-hp-II-sm} (for {\bf Case}
(i)) and Table \ref{M-hp-II-nsm} (for {\bf Case} (ii)).

\vskip 0.1in
\begin{center}
\tabcaption{}\label{M-hp-II-sm}
 Iteration counts of PCG with the preconditioner $B_{II}$ (hexahedral elements): the coefficients have no jump
    \vskip 0.2in
\begin{tabular}{|c|c|c|c|c|c|}
\hline
$ m \backslash n$ & ~~4~~ &~~6~~  &~~8~~   &~~10~~ \\\hline
 8& 23 &22 & 21 &21    \\\hline
 16&23 &23 & 22 &21     \\\hline
 24&25 &25 &24  &23     \\\hline
 32&26 &26 &25  &24    \\\hline
\end{tabular}
\end{center}
 \vskip 0.1in

\vskip 0.1in
\begin{center}
\tabcaption{}\label{M-hp-II-nsm}
 Iteration counts of PCG with the preconditioner $B_{II}$ (hexahedral elements): the coefficients have large jumps
\vskip 0.2in
\begin{tabular}{|c|c|c|c|c|c|c|c|c|}\hline
$\mbox{}$& \multicolumn{4}{c|}{Choice~(1)~of~$D$}&
\multicolumn{4}{c|}{Choice~(2)~of~$D$}\\\hline $\mbox{}$&
\multicolumn{2}{c|}{$ \beta_0 =\alpha_0 =10^{-5}  $
}&\multicolumn{2}{c|}{$ \beta_0 =\alpha_0=10^{5} $}&
\multicolumn{2}{c|}{$ \beta_0 =\alpha_0=10^{-5}  $
}&\multicolumn{2}{c|}{$ \beta_0 =\alpha_0=10^{5} $}\\\hline
$~m ~ \backslash n$& ~~4~~&~~8~~&~~4~~&~~8~~&
~~4~~&~~8~~&~~4~~&~~8~~\\\hline
 8&  18   &21  &27  &25    &18    &21  &27  & 29     \\\hline
 16& 19   &21  &29  &26    &19    &21  &29  & 28    \\\hline
 24& 21   &23  &31  &28    &21    &24  &31  & 31  \\\hline
 32& 21   &25  &33  &30    &21    &25  &33  & 34  \\\hline
\end{tabular}
\end{center}
 \vskip 0.1in

From above table, we observe that the  convergence rate of PCG method with the preconditioner $B_{II}$ is quasi-optimal, even if the coefficients
have large jumps.

Next we consider the case of the tetrahedral elements. We use the same way to define a coarse space as in Subsection 6.2.1 for this case.
We list the iteration counts of PCG method with the preconditioner $B_{II}$ in Table \ref{M-tp-II-sm} (for {\bf Case (i)}) and Table \ref{M-tp-II-nsm}
 (for {\bf Case (ii)}).

\vskip 0.1in
\begin{center}
 \tabcaption{}\label{M-tp-II-sm}
 Iteration counts of PCG with the preconditioner $B_{II}$ (tetrahedral elements): the coefficients have no jump
    \vskip 0.2in
\begin{tabular}{|c|c|c|c|c|c|}
\hline
$ m \backslash n$ & ~~4~~ &~~6~~  &~~8~~   &~~10~~ \\\hline
 8& 21 &21 & 20 &20    \\\hline
 16&26 &25 & 24 &24     \\\hline
 24&26 &25 & 24 &24       \\\hline
 32&28 &28 &27 &26    \\\hline
\end{tabular}
\end{center}
 \vskip 0.1in

\vskip 0.1in
\begin{center}
\tabcaption{}\label{M-tp-II-nsm}
 Iteration counts of PCG with the preconditioner $B_{II}$ (tetrahedral elements): the coefficients have large jumps
\vskip 0.2in
\begin{tabular}{|c|c|c|c|c|c|c|c|c|}\hline
$\mbox{}$& \multicolumn{4}{c|}{Choice~(1)~of~$D$}&
\multicolumn{4}{c|}{Choice~(2)~of~$D$}\\\hline $\mbox{}$&
\multicolumn{2}{c|}{$ \beta_0 =\alpha_0 =10^{-5}  $
}&\multicolumn{2}{c|}{$ \beta_0 =\alpha_0=10^{5} $}&
\multicolumn{2}{c|}{$ \beta_0 =\alpha_0=10^{-5}  $
}&\multicolumn{2}{c|}{$ \beta_0 =\alpha_0=10^{5} $}\\\hline
$ m\backslash n$& ~~4~~&~~8~~&~~4~~&~~8~~&~~4~~&~~8~~&~~4~~&~~8~~\\\hline
 8& 18    &20  &27  &24    &18    &20  &26   & 26     \\\hline
 16& 21   &24  &32  &29    &21    &24  &32  & 31    \\\hline
 24& 21   &23  &32  &28    &21    &23  &31  & 31        \\\hline
 32&23    &27  &35  &32    &23    &27  &35  & 35  \\\hline
\end{tabular}
\end{center}
 \vskip 0.1in

It can be seen from the above two tables that the iteration counts of the PCG method with the new preconditioner only slowly grow when $m=d/h$ increases, but not depend on $n=1/d$.

\subsection{On the proposed coarsening technique}
It can be seen, from the results in Subsection 6.1 and Subsection 6.2, that the preconditioner $B_{II}$ has almost the same convergence rate as
the preconditioner $B_{I}$. A key ingredient in the preconditioner $B_{II}$ is the proposed coarsening technique.
In this subsection, we give some numerical result to illustrate the efficiency of the coarsening technique and further explain that the preconditioner $B_{II}$ is indeed very cheap.

We first consider a typical cuboid domain $G=[0,2]\times [0,1] \times [0,1]$ to investigate the approximate effect of the coarsening technique.
Let $G$ be divided into the union of two cube $G_1$ and $G_2$, with $G_1=[0, 1]^3$ and $G_2=[1,2]\times [0,1] \times [0,1]$, and set $\ff=\partial G_1\cap\partial G_2$.
We will compare the accuracy of the solutions of the two systems (\ref{4.new3}) (with the original partition ${\mathcal T}_h$) and (\ref{4.new4}) (with the coarsening partition $\tilde{\mathcal T}^{\ff}_{\tilde{h}}$),
where $\Omega_{l}^{\ff} = G_l~(l = 1, 2)$. To this end, we need to calculate the discrete $l^2$ relative error on $\ff$, which is defined by
$$ err. = \frac{||\tilde{\chi}_{\ff}- \chi_{\ff}||_{l^2}}{||\chi_{\ff}||_{l^2}}.$$

In  Table \ref{l2-re-error}, we list the results for the linear elasticity problem and Maxwell's equations, with constant coefficients.
   \vskip 0.1in
\begin{center}
  \tabcaption{}\label{l2-re-error}
  The $l^2$ relative error restricted on the interface $\ff$ for coarsening
\vskip 0.2in
\begin{tabular}{|c|c|c|}
\hline
$\mbox{}$ &  \multicolumn{2}{c|}{$err.$}\\\hline
$h$      & Linear elasticity problem & Maxwell's equations
           \\\hline

 1/8&    0.0407     &0.0582       \\\hline
 1/16&   0.0240     &0.0317      \\\hline
  1/24&   0.0225      &0.0284      \\\hline
 1/32&  0.0233     &0.0254       \\\hline
\end{tabular}
\end{center}
 \vskip 0.1in

From this table, we can see that the solution of the auxiliary system (\ref{4.new4}) indeed is a rough approximation of the solution of the interface system (\ref{4.new3}),
which can explain why the preconditioner $B_{II}$ is effective, as confirmed in Subsection 6.1.2 and Subsection 6.2.2.

Next, we illustrate the local solvers in $B_{II}$ indeed is very cheap. For simplicity, we just select one face $\ff$ and one interior vertex $\vv \in \mathcal{N}_d$ to test our coarsening technique,
where the face $\ff$ is shared by $\Omega_i$ and $\Omega_j$. Let $n_c$ and $n_f$ denote the dofs corresponding to the coarsening partition and the original
fine partition on $\Omega_{\ff}$ (or $\Omega_{\vv}$), respectively.

In  Table \ref{coarse-for-II-nodal} and  Table \ref{coarse-for-II-edge}, we list the local dofs $n_c$ and $n_f$ for linear elasticity problem and Maxwell's equations, respectively.

\vskip 0.1in
\begin{center}
  \tabcaption{}\label{coarse-for-II-nodal}
 The dofs of local problems solved in Step 3 and Step 4 of {\bf Algorithm 4.1}: linear elasticity problem (with vector-valued nodal basis functions)
\vskip 0.2in
\begin{tabular}{|c|c|c|c|c|c|c|}
\hline
$\mbox{}$ & \multicolumn{3}{c|}{$\Omega_{\ff}$} &
            \multicolumn{3}{c|}{$\Omega_{\vv}^{half}$ }\\\hline
$d/h$      & $n_c$ (coarse) & $n_f$ (fine) &$n_c/n_f$
           & $n_c$ (coarse) & $n_f$ (fine) &$n_c/n_f$ \\\hline

 8&    173*3    &735*3   &0.24  &65*3     &537*3   &0.12     \\\hline
 16&   809*3   &6975*3   &0.12  &306*3    &4145*3  &0.07     \\\hline
 24&   1951*3   &24863*3 &0.08  &749*3    &13897*3 &0.05     \\\hline
 32&   3585*3  &60543*3  &0.06 &1245*3    &33729*3 &0.04      \\\hline
\end{tabular}
\end{center}
 \vskip 0.1in

\vskip 0.1in
\begin{center}
 \tabcaption{}\label{coarse-for-II-edge}
  The dofs of local problems solved in Step 3 and Step 4 of {\bf Algorithm 4.1}: Maxwell's equations (with Nedelec edge basis functions)
\vskip 0.2in
\begin{tabular}{|c|c|c|c|c|c|c|}
\hline
$\mbox{}$ & \multicolumn{3}{c|}{$\Omega_{\ff}$} &
            \multicolumn{3}{c|}{$\Omega_{\vv}^{half}$ }\\\hline
$d/h$      & $n_c$  (coarse) & $n_f$ (fine) &$n_c/n_f$
           & $n_c$  (coarse) & $n_f$ (fine) &$n_c/n_f$ \\\hline

 8&    1673    &6240   &0.38  &1910    &5554      &0.34     \\\hline
 16&   6642    &53568  &0.18  &5390    &35658     &0.15     \\\hline
 24&   15157   &184992 &0.08  &11076   &111842    &0.10             \\\hline
 32&   27132   &443520 &0.06  &16575   &255610    &0.06
      \\\hline

\end{tabular}
\end{center}
 \vskip 0.1in

It can be seen from these results that the dofs of local problems in Step 3 and Step 4 of {\bf Algorithm 4.1} are much smaller than that
of the local problems associated with the original fine grids. Moreover, the smaller the value $d/h$ is, the better the coarsening effect is.
This means that the preconditioner $B_{II}$ is very cheap.

Finally, we illustrate why this coarsening technique has not been applied directly to the preconditioner $B_I$ (see {\bf Remark 4.2}).

In Table \ref {coarse-for-I}, we list the dofs of the problems in Step 3 of {\bf Algorithm 3.1}, which associated with the fine partition and the coarsening partition on $\Omega_{\vv}^{half}$, respectively.

\vskip 0.1in
\begin{center}
 \tabcaption{}\label{coarse-for-I}
   The dofs of the local problems solved in Step 3 of {\bf Algorithm 3.1}: coarsening or not
\vskip 0.2in
\begin{tabular}{|c|c|c|c|c|c|c|}
\hline
$\mbox{}$ & \multicolumn{3}{c|}{Elasticity problems} &
            \multicolumn{3}{c|}{Maxwell's problems }\\\hline
$d/h$      & $n_c$  (coarse) & $n_f$ (fine) &$n_c/n_f$
           & $n_c$ (coarse) & $n_f$ (fine) &$n_c/n_f$ \\\hline
 8&    289*3     &792*3    &0.40  &3048    &6130   &0.50     \\\hline
 16&   1591*3    &4913*3   &0.32  &14704   &37962  &0.39     \\\hline
 24&   3964*3    &15625*3  &0.25  &32115   &117026  &0.27     \\\hline
 32&   7256*3    &35937*3  &0.20  &60118   &264826  &0.23
      \\\hline

\end{tabular}
\end{center}
 \vskip 0.1in

From  Table \ref{coarse-for-I}, we can see that the coarsening effect is not ideal for this situation.

\section{Conclusion}
 In this paper, we have
constructed two substructuring preconditioners with the simplest coarse space for
general elliptic-type problems in three dimensions. In particular, we design two kinds of
new local interface solvers, which are easy to implement and do not depend on the considered models.
The proposed preconditioners can absorb some advantages of the non-overlapping DDMs
and the overlapping DDMs. Especially, in the second
preconditioner we propose a coarsening technique to solve local interface problems.
As expected, the utilization of coarsen grids can significantly
decrease the cost of calculation, but does not destroy the convergence
rate of the PCG method.
We have given some numerical results to show that the proposed preconditionners are effective uniformly for
the linear elasticity problem and Maxwell's equations in three dimensions.


\end{document}